\documentclass[reqno,10pt]{article}
\usepackage{amsmath}
\usepackage{graphicx, color, enumerate}
\usepackage{amscd}
\usepackage{amsmath}
\usepackage{amsfonts}
\usepackage{amssymb}
\usepackage{mathrsfs}
\usepackage{latexsym}

\newcommand{\be}{\begin{equation}}
\newcommand{\ee}{\end{equation}}

\newtheorem{corollary}{Corollary}[section]
\newcommand{\R}{\mathbb{R}}
 \newcommand{\Rn}{\mathbb{R}^n}

\newcommand{\E}{\mathbb{E}}

\newcommand{\dyle}{\displaystyle}

\renewcommand{\a }{\alpha }
\renewcommand{\b }{\beta }

\newcommand{\D }{\Delta }
\newcommand{\e }{\varepsilon }

\renewcommand{\l }{\lambda }

\newcommand{\m }{\mu }
\newcommand{\n }{\nabla }

\newcommand{\s }{\sigma }

\newcommand{\cN}{{\cal N}}

\newcommand{\Wan}{W^{1,2}(\Rn)}

\newcommand{\intn}{\int_{\Rn}}

\newcommand{\wt}{\widetilde}

\newcommand{\bE}{{\bf E}}
\newcommand{\bu}{{\bf u}}
\newcommand{\bv}{{\bf v}}
\newcommand{\bz}{{\bf z}}

\newcommand{\bo}{{\bf 0}}

\numberwithin{equation}{section}

\newtheorem{Theorem}{Theorem}

\newtheorem{Definition}[Theorem]{Definition}
\newtheorem{remark}[Theorem]{Remark}
\newtheorem{remarks}[Theorem]{Remarks}
\newtheorem{example}[Theorem]{Example}
\newtheorem{examples}[Theorem]{Examples}
\newenvironment{Remark}{\begin{remark} \rm}{\rule{2mm}{2mm}\end{remark}}

\numberwithin{Theorem}{section}

\begin{document}
\title{
{\bf\Large  Positive solutions to some systems of coupled nonlinear Schr\"odinger equations}}
\author{{\bf\large Eduardo Colorado}\hspace{2mm}
{\bf\large}\vspace{1mm}\\
{\it\small Departamento de Matem\'aticas, Universidad Carlos III de Madrid}\\
{\it\small Avda. Universidad 30, 28911 Legan\'es (Madrid), Spain,}\\
{\it\small and}\\
{\it\small Instituto de Ciencias Matem\'aticas, ICMAT (CSIC-UAM-UC3M-UCM)}\\
{\it\small C/Nicol\'as Cabrera 15, 28049 Madrid, Spain.}\\
{\it\small e-mail: eduardo.colorado@uc3m.es, eduardo.colorado@icmat.es}}

\maketitle

\begin{center}
{\bf\small Abstract}

\vspace{3mm}
\hspace{.05in}\parbox{4.5in}
{{\small We study the existence of nontrivial bound state  solutions to the following system of coupled nonlinear time-independent Schr\"odinger equations
$$
- \D u_j+ \l_j u_j =\m_j u_j^3+ \sum_{k=1;k\neq j}^N\b_{jk}  u_ju_k^2,\quad u_j\in \Wan;\: j=1,\ldots,N
$$
where  $n=1,\,2,\, 3; \,\l_j,\, \m_j>0$ for $j=1,\ldots ,N$, the coupling parameters $\b_{jk}=\b_{kj}\in \mathbb{R}$ for $j,k=1,\ldots,N$, $j\neq k$.
Precisely, we  prove the existence of  nonnegative bound state solutions for suitable conditions on the coupling factors.
Additionally, with more restrictive conditions on the coupled parameters, we show that the bound states founded are positive.
}}
\end{center}

\vglue 0.5cm

\noindent{\it \footnotesize 2010 Mathematics Subject Classification}. {\scriptsize  35Q55, 35B38, 35J50, 47A55}.\\
{\it \footnotesize Key words}. {\scriptsize Nonlinear Schr\"odinger Equations,  Bound States, Critical Point Theory, Perturbation Theory}

\tableofcontents
\section{Introduction}
Problems on coupled Nonlinear Schr\"odinger Equations (NLSE for short) have been widely investigated in the last years. They arise
naturally in nonlinear Optics, and in the Hartree-Fock theory for Bose-Einstein condensates, among other physical problems.
For example, a planar light beam propagating in the $z$ direction in a non-linear medium can be described by a NLSE  of the form
$$
{\rm i}\, \bE_z+ \bE_{xx}+ \theta |\bE|^2 \bE=0,
$$
where ${\rm i}$, $\bE(x,z)$  denote the imaginary unit and  the complex envelope of an Electric field respectively, and $\theta>0$ is a normalization constant, corresponding to the fact that the medium is self-focusing.

Here we consider the following more general  system of $N$-coupled   NLSE
\be\label{eq:main-sys-time}
\left\{\begin{array}{ll}
&\dyle -{\rm i}\, \frac{\partial}{\partial t} E_j- \D E_j = \mu_j |E_j|^2E_j+\sum_{k\neq j;k=1}^N\b_{jk}  |E_k|^2E_j,\quad x\in\Rn,\quad t>0\\ & \\
& E_j=E_j(x,t)\in \mathbb{C},\quad E_j(x,t)\to 0  \mbox{ as }|x|\to \infty
\end{array}\right.
\ee
for $j=1,\ldots,N$, the coupled parameters $\b_{jk}=\b_{kj}\in \R$ for $j,k=1,\ldots,N$, $j\neq k$; $\l_j>0$ and  $\mu_j>0$ is for self-focusing in the $j$-th component of the beam. The solution $E_j$ represents the $j$-th component of the beam. The coupling constant $\b_{jk}$ means the interaction between the $j$-th and $k$-th component of the beam. If they are positive, the interaction is attractive while if they are negative we have a repulsive interaction. The mixed case involves different sign on the coupling parameters $\b_{jk}$. Here, we study  the attractive and mixed interactions.

When one looks for solitary wave solutions of the form $E_j(x,t)=e^{i\l_jt}u_j(x)$, with $\l_j>0$, $u_j$ are the real valued functions called standing wave solutions which solve the following system of $N$-coupled nonlinear time-independent Sch\-r\"o\-dinger equations
\be\label{eq:main-sys}
- \D u_j+ \l_j u_j = \mu_j u_j^3+\sum_{k\neq j;k=1}^N\b_{jk}  u_k^2u_j, \qquad u(x)\to 0 \mbox{ as }|x|\to \infty,
\ee
for $j=1,\ldots,N$. Solutions $u_j$  belong to the Sobolev space $\Wan$ ($n=1,\,2,\,3$) for all $j=1,\ldots,N$.

In the last years there has been a very active investigation on coupled systems of NLSE, from the physical point of view; see for instance \cite{hioe,me} or the book \cite{akanbook} and the references therein, as well as many works of mathematicians, from a more theoretical point of view,  dealing with the existence, multiplicity, uniqueness and qualitative properties of bound and ground state solutions; see the earlier works \cite{ac1,ac2,acr,bt,linwei,mmp,pomp,sirakov,ww}, and the more recent list (far from complete) \cite{chen-zou,dancer-wei,itanaka,liu-wang,soave,wei-yao}. This kind of systems have also been recently studied in the framework of the fractional and bi-harmonic Schr\"odinger operators in \cite{col} and \cite{acv} respectively.

Most of the existing works have been developed for systems with two equations, i.e., $N=2$. For a general $N$-system, in \cite{linwei} Lin-Wei  studied the attractive interaction case ($\b_{ij}>0$) in the dimensional case $n=1,\,2,\,3$ and small coupling factors, i.e., when the ground state could coincide with a semitrivial solution, which corresponds to the ground state of a single NLS equation in one component and zero the others. More recently, in \cite{liu-wang}  Liu-Wang gave a sufficient condition on large coupling coefficients for the existence of a nontrivial ground state solution in a system of $N$ ($N$-system for short) NLSE for the dimensional case $n=2,\, 3.$

The existence of positive bound states was firstly studied by Ambrosetti-C. in \cite{ac1,ac2}, where for a $2$-system in the dimensional case $n=2,3$  was proved, among other results, the existence of positive bound states for positive non-large coupling factors; see also \cite{col} for $N=2,\,n=1$. Even more, that existence result was also proved in \cite{ac2} for $N\ge 3$ in the attractive interaction case, provided the coupling factors are sufficiently small. Later in \cite{itanaka} with $N=2$, Ikoma-Tanaka shown that bound states obtained by Ambrosetti-C. in \cite{ac1,ac2} are indeed least energy (positive) solutions in some range of the small coupling factor. That result was recently improved by Chen-Zou in \cite{chen-zou} for a better suitable range of the parameters.

Here, we improve the only result (up to our knowledge) of the existence of bound state solutions for a system of coupled  NLSE such us \eqref{eq:main-sys} with $N\ge 3$. Precisely, we demonstrate the existence of new positive and non-trivial bound state solutions to \eqref{eq:main-sys} in dimensions $n=1,2,3$ different from all previously studied. We want to emphasize that the interaction  between the components of the beam  is analyzed here in the attractive case  and the mixed one.

The paper is organized as follows. In Section \ref{sec:not-prel} we introduce the functional setting, the notation and give  some definitions. Section \ref{sec:pr} contains some previous known results, while the last one, Section \ref{sec:bs} is devoted to prove the main results of the work.

\section{Functional setting and notation}\label{sec:not-prel}
In this section, as we just mentioned above, we will establish the corresponding functional framework, the notation and give some definitions. The functional setting comes from the variational structure of \eqref{eq:main-sys}, which can be derived from
$$
- \D u_j + \l_j u_j =   f_{u_j}(u_1,\ldots,u_n)
$$
where
$$
f(u_1,\ldots,u_n)=\dyle\sum_{j,k=1}^N\a_{jk}u_j^2u_k^2,
$$
with $\a_{jj}=\frac 14\mu_j$ for $j=1,\ldots,N$; the coupling parameters $\a_{jk}=\frac 14\b_{jk}$ are symmetric,  $\b_{jk}=\b_{kj}$, both for $k,j=1,\ldots,N$ with $k\neq j$.
We will work on the Sobolev  functional space $\Wan$, i.e., we will assume that $u_j\in \Wan$ for $j=1,\ldots,N$. We remember that the Sobolev space $E=\Wan$ can be defined as the completion of $\mathcal{C}_0^1(\Rn)$ with the norm given by
$$
\|u\|_{E}=\left(\intn ( |\n u|^2+ u^2)\, dx\right)^{\frac 12},
$$
which comes from the scalar product
$$
\langle u|v\rangle_E=\intn (\n u\cdot \n v+ uv)\, dx.
$$
We will denote the following equivalent norms in $E$,
\be\label{norma}
\|u\|=\|u\|_j=\left(\intn (|\n u|^2+\l_j u^2)\, dx\right)^{\frac 12},
\ee
with the corresponding scalar products
\be\label{p-escalar}
( u|v)=( u|v)_j=\intn (\n u\cdot \n v+ \l_j uv)\, dx; \quad j=1,\ldots,N,
\ee
dropping the sub-index $j$ in \eqref{norma}-\eqref{p-escalar} for short.
In this manner, solutions of \eqref{eq:main-sys} are the critical points $\bu=(u_1,\ldots,u_N)\in \E=E\times\stackrel{(N)}{\cdots}\times E$ of the  corresponding energy functional defined by
$$
\Phi (\bu)=\frac 12\|\bu\|^2- F(\bu),
$$
where
$$
F(\bu)= \intn f(\bu)\, dx\quad {and}\quad
\|\bu\|^2=\sum_{j=1}^N \|u_j\|^2.
$$
Also we define
$$
I_j(u)=\frac 12\intn (|\n u|^2+\l_j u^2)\, dx-\frac 14 \m_j\intn u_j^4\, dx\quad\mbox{for }j=1,\ldots,N
$$
With respect to the  coupling factors, we assume $m\in \mathbb{N}$, moreover $2\le 2m\le N$, and define $\b_k=\b_{2k-1,2k}$ and $\l_{2k-1}=\l_{2k}$ for $1\le k\le m$. If $2m<N$,  we define
 $\b_{m+\ell}^{(2k)}=\b_{2m+\ell,2k}$, $\b_{m+\ell}^{(2k-1)}=\b_{2m+\ell,2k-1}$  for any $1\le \ell \le N-2m$ and $k\neq 2m+\ell$. We also suppose that $\b_{m+\ell}^{(2k)}=\e\wt{\b}^{(2k)}_{m+\ell}$, $\b_{m+\ell}^{(2k-1)}=\e\wt{\b}^{(2k-1)}_{m+\ell}$.

 Let us denote $\Phi_\e=\Phi$ to emphasize its dependence on $\e$, then we can rewrite the functional $\Phi_\e$ as
$$
\Phi_{\e}(\bu)=\Phi_0-\e \wt{F} (\bu)
$$
where
$$
\Phi_0(\bu)=\sum_{j=1}^N I_j(u_j)-\frac 12\sum_{k=1}^m\b_{k}   \intn u_{2k-1}^2u_{2k}^2\,dx
$$
and
$$
\wt{F}(\bu)=\frac 12 \sum_{\ell=1}^{N-2m} \sum_{k=1}^m\intn \left(\wt{\b}_{k+\ell}^{(2k-1)}u_{2k-1}^2+\wt{\b}_{k+\ell}^{(2k)}u_{2k}^2\right) u_{2m+\ell}^2\,dx.
$$
For a function of two components we also set
\be\label{eq:Phi-k}
\Phi_k(u,v)=I_{2k-1}+I_{2k}-\frac 12 \b_k\intn u^2v^2\,dx\qquad 1\le k\le m.
\ee
We denote $\bo=(0,\ldots,0)$, then we say that a vector function (or a function to simplify) $\bu\in\E$ is positive (non-negative), namely $\bu>\bo$ ($\bu\ge \bo$) if every component of $\bu$ is positive (non-negative), i.e., $u_j>0$ ($u_j\ge 0$) for every $j=1,\ldots,N.$

\begin{Definition}
A bound state $\bu\in \E$ of (\ref{eq:main-sys}) is a critical point
of $\Phi$. A nontrivial bound state $\bu=(u_1,\ldots,u_N)$ is a bound state with all the components $u_j\not\equiv 0$, $j=1,\ldots,N$. Moreover, a positive bound state   $\bu\ge \bo$  such that its energy is minimal among {\it all the non-trivial bound states}, namely
\begin{equation}\label{eq:gr}
\Phi(\bu)=\min\{\Phi(\bv): \bv\in \E\setminus\{\bo\},\; \Phi'(\bv)=0\},
\end{equation}
is called a {\it ground state}  of \eqref{eq:main-sys}.
\end{Definition}

Note that a ground state can have some of the components equal zero, thus in some works the definition of a ground state corresponds to a positive ground state here.

We will focus on the existence of (non-trivial) non-negative and also on positive bound states.

Let us  denote $2^*=\frac{2n}{n-2}$ if $n\ge 3$ and $2^*=\infty$ if $n=1,\,2$, and remark that the functional make sense because $E\hookrightarrow L^4(\Rn)$ for $n=1,2,3$. Concerning the Palais-Smale $(PS)$ condition, one cannot expects a compact embedding neither in $E$ nor in the restriction on the even functions of $E$ into $L^q(\R)$ for any $1<q<\infty$.  If $n=2,\, 3$, the PS condition is easy to obtain because the space of radially symmetric functions of the Sobolev space  $E$, namely $H=E_{rad}$,  is compactly embedded into $L^{q}(\mathbb{R}^n)$ for $q<2^*$ (in particular for $q=4$), see for instance \cite{lions}.
Nevertheless, as we will see, the lack of compactness will not be a problem provided we will work on $H$  for all $n=1,2,3$. The reason is that we can obtain strong convergence in
$\mathbb{H}=H\times\stackrel{(N)}{\cdots}\times H$ thanks to the Local Inversion Theorem by a perturbation argument involving some non-degenerate critical points.
\section{Previous results}\label{sec:pr}
In order to prove the main result of the paper we introduce some known results we are going to use.
\begin{Remark}\label{rem:1}
Let $U$ be the unique positive solution in $H$ of the equation $- \D u+u=u^3$ (see the paper of Kwong \cite{kw}) then the scaled function
\be\label{eq:single}
U_j(x)=\sqrt{\frac{\l_j}{\mu_j}}\,U(\sqrt{\l_j}\,x),\quad j=1,\ldots,N
\ee
is the unique positive solution  in $H$ of  $- \D u+\l_ju=\m_j u^3$ for $j=1,\ldots,N$. We also recall that $U_j$ is the ground state, as well as a non-degenerate critical point of the corresponding functional $I_j$, $j=1,\ldots, N$.
\end{Remark}
Let us consider the following $2$-system of NLSE
\be\label{eq:2-sistema}
(S_2)\equiv\left\{\begin{array}{ll}
- \D u + \l u &= \mu_1 u^3+\b  v^2 u\\
- \D v + \l u &= \mu_2 v^3+\b u^2 v.
\end{array} \right.
\ee
Throughout this subsection we will maintain the notation of Section \ref{sec:not-prel} for a $2$-system  or a $3$-system of  NLSE with the natural meaning.

Let us  define
$$
 a_{k\b} = \sqrt{\frac{\mu_k(\mu_j-\b)}{\mu_k\mu_j
-\b^2}},\quad k\not= j,\;k=1,2.
$$
If $\l_1=\l_2=\l$,  we have the following explicit and positive solution of $(S_2)$, given by
\be\label{eq:2-ground}
(u^0,v^0)=(a_{1\b}U_1,a_{2\b}U_2)
\ee
with $U_1,\, U_2$ defined  in \eqref{eq:single}, see also \cite[Remarks 5.8-(b)]{ac2}.
\begin{Remark}\label{rem:non-degenerate}
  In  \cite[Lemma 2.2 and Theorem 3.1]{dancer-wei} Dancer-Wei proved that $(u^0,v^0)$ is a non-degenerate critical point to the corresponding energy functional of $(S_2)$.
\end{Remark}
Even more, Wei-Yao  proved in \cite{wei-yao} the following property of the critical point $(u^0,v^0)$.

\

\noindent
{\bf Theorem A.} {\it Assume that $\b>\max\{\m_1,\m_2\}$}, then $(u^0,v^0)$ is the unique positive solution of $(S_2).$
\begin{Remark}
According to \cite[Remark5.8-(b)]{ac2} and Theorem A one has that $(u^0,v^0)$ is a ground state of $(S_2)$, i.e.,
$$
\Phi(u^0,v^0)=\min\{\Phi(\bv): \bv\in \E\setminus\{\bo\},\; \Phi'(\bv)=0\},
$$
which also corresponds to
$$
\Phi(u^0,v^0)=\min_{\bu\in\cN}\Phi(\bu)
$$
where $\cN$ is the corresponding Nehari manifold, defined as
$$
\cN =\{ \bu\in \mathbb{E}\setminus\{{\bf 0}\}: (\Phi'(\bu)|\bu)=0\},
$$
that is  a natural constraint for $\Phi$; see \cite{ac2} for more details. Furthermore,
\be\label{eq:inf1}
\begin{array}{rcl}
 \|(u^0,v^0)\|& = &\dyle \inf_{\bu\in \mathbb{H}\setminus\{\bo\}} \frac{\|(u_1,u_2)\|^2}{\dyle\left(\m_1 \intn u_1^4\,dx+\m_2 \intn u_2^4\,dx+2\b\intn u_1^2u_2^2\,dx\right)^{\frac12}}\\ & &  \\
 & = & \dyle\inf_{\bu\in \mathbb{E}\setminus\{\bo\}} \frac{\|(u_1,u_2)\|^2}{\dyle\left(\m_1 \intn u_1^4\,dx+\m_2 \intn u_2^4\,dx+2\b\intn u_1^2u_2^2\,dx\right)^{\frac12}}.
\end{array}\ee
\end{Remark}
\section{Existence of bound states}\label{sec:bs}
In order to prove the existence  of a bound state solution to System \eqref{eq:main-sys}, first we demonstrate  a reduced result to both, clarify and also simplify the proof of the main result in Theorem \ref{th:ec2}. To do so we consider the following $3$-system of NLSE
$$
(S_3)\equiv \left \{
\begin{array}{ll}
- \D u_1+ \l_1 u_1 &= \mu_1 u_1^3+\b_{12}  u_1u_2^2+ \b_{13} u_1u_3^2\\
- \D u_2 + \l_2 u_2 &= \mu_2 u_2^3+\b_{21} u_2 u_1^2 + \b_{23}  u_2u_3^2\\
- \D u_3 + \l_3 u_3 &= \mu_3 u_3^3+\b_{31} u_3u_1^2 + \b_{32}u_3u_2^2
\end{array} \right.
$$
with $\l_j,\m_j>0$, $u_j\in H$, $j=1,\,2,\,3$ and $\b_{jk}=\b_{kj}$ for all $j,k=1,\,2,\, 3$; $j\neq k$.
\begin{Remark}\label{rem:sobolev-cte}
Let
  $$
\s_\l^2\equiv \inf_{\varphi\in E\setminus \{ 0\}}\frac{\|
\varphi\|_\l^2}{\| \varphi\|^2_{L^4}}=\inf_{\varphi\in E,\, \|
\varphi\|_{L^4}=1}\| \varphi\|_\l^2,
$$
denotes  the  best Sobolev constant in the embedding
of $(E,\|\cdot\|_{\l})$ into  $L^4(\Rn)$, where
$$
\| u\|^2_\l=\intn (|\n u|^2+\l u^2)\, dx.
$$
Easily one obtains that  $\s_\l$ is attained  at
$$
v_\l(x)=\s_\l^{-1}\sqrt{\l} U(\sqrt{\l}x),
$$
and one has
$$
\s_\l^4=\l^2\intn U^4(\sqrt{\l}x)dx=\l^{2-\frac n2}\intn U^4(x)dx.
$$
\end{Remark}
Let us set $\bu_0=(u^0,v^0,U_3)$, where $(u^0,v^0)$, $U_3$ are given by \eqref{eq:2-ground}, \eqref{eq:single} respectively. Then we have the following.
\begin{Theorem}\label{th:ec1}
If $\l_1=\l_2=\l$, $\b_{12}>\max\{ \m_1,\m_2\}$, 
and $\b_{j3}=\e\wt{\b}_{j3}$, $j=1,2$, there exists $\e_0>0$ such that for $0<\e<\e_0$,  $(S_3)$ has a radial bound state $\bu_{\e}\ge \bo$ with  $\bu_{\e}\to \bu_0$ as $\e\to 0$. Furthermore, if both  $\b_{13}, \,\b_{23}\ge 0$ then $\bu_{\e}>\bo$.
\end{Theorem}
Notice that this result corresponds to $N=3$, $m=1$,  $\ell=1$.
\begin{Remark}\label{rem:2}
Note that hypothesis $\b_{12}>\max\{ \m_1,\m_2\}$ comes from Theorem A, on the contrary $(S_2)$ have no positive solutions. Multiplying the first (third) equation of $(S_3)$ by $u_3$ ($u_1$) and integrating on $\Rn$ we get
\be\label{eq:h1}
(\m_1-\b_{13})\intn \!\! u_1^3u_3\, dx+(\b_{12}-\b_{23})\intn\!\! u_1u_2^2u_3\, dx +(\b_{13}-\m_3)\intn\!\! u_1u_3^3\, dx=0.
\ee
Multiplying the second (third) equation of $(S_3)$ by $u_3$ ($u_2$) and integrating again on $\Rn$ we find
\be\label{eq:h2}
(\m_2-\b_{23})\intn\!\! u_2^3u_3\, dx+(\b_{12}-\b_{13})\intn\!\! u_1^2u_2u_3\, dx +(\b_{23}-\m_3)\intn\!\! u_2u_3^3\, dx=0.
\ee
As a consequence, if for instance $\mu_3<\b_{13},$ thanks to \eqref{eq:h1}, system $(S_3)$ has no positive solutions. Similarly, if $\m_3<\b_{23}$, due to   \eqref{eq:h1}, system  $(S_3)$ has no positive solutions. But $\m_3>\max\{\b_{13},\,\b_{23}\}$  provided $\e_0$ is sufficiently small.
\end{Remark}
\noindent {\it Proof of Theorem \ref{th:ec1}}.
Taking  $\e_0>0$ sufficiently small, we can suppose that $\m_3>\max\{\b_{13},\, \b_{23}\}$ and the compatibility conditions to have non-negative solutions given in \eqref{eq:h1}, \eqref{eq:h2} are satisfied. According to the notation of Section \ref{sec:not-prel} with $N=3$, $m=1$,  $\ell=1$, we  have
$$
\Phi_{\e}(\bu)=\Phi_0-\e \wt{F} (\bu)
$$
where
$$
\Phi_0(\bu)=\sum_{j=1}^3 I_j(u_j)-\frac 12\b \intn u_1^2u_2^2\,dx
$$
and
$$
\wt{F}(\bu)=\frac 12 \intn (\wt{\b}_{13}u_1^2u_3^2+\wt{\b}_{23}u_2^2u_3^2)\,dx.
$$
Let us consider the critical point $\bu_0$ of the unperturbed functional $\Phi_0$. We note that $U_3$ is a non-degenerate critical point of $I_3$ on $H$; see \cite{kw} and Remark \ref{rem:1}. Also, as we pointed out in Remark \ref{rem:non-degenerate}, $(u^0,v^0)$  is a non-degenerate critical point to the corresponding energy functional (of $(S_2)$ with $\l_1=\l_2=\l$ and $\b_{12}=\b$) acting on $\mathbb{H}$. As an immediate consequence, $\bu_0$ is a non-degenerate critical point of $\Phi_0$ on $\mathbb{H}$, then a straightforward application of the Local Inversion Theorem yields the existence of a critical point $\bu_\e$ of $\Phi_\e$ for any $0<\e<\e_0$ with $\e_0$ sufficiently small; see  \cite{a-m} for more details. Moreover, $\bu_\e\to \bu_0$ on $\mathbb{H}$ as $\e\to 0$. To complete the proof it remains to show that $\bu_\e\ge \bo$, and furthermore, in case $\b_{13},\,\b_{23}\ge 0$ then $\bu_\e>\bo$. We follow an argument of \cite{cing} with suitable modifications by separating the positive and negative parts and using energy type estimates.

Let us  denote the positive part
 $\bu_\e^+=(u_{1\e}^+,u_{2\e}^+,u_{3\e}^+)$ and  the negative part $\bu_\e^-=(u_{1\e}^-,u_{2\e}^-,u_{3\e}^-)$.
Notice that $(u^0,v^0)$ satisfies the identity \eqref{eq:inf1} and $U_3$ satisfies the following identity
\be\label{eq:inf2}
 \|U_3\|=\inf_{u\in H\setminus\{0\}} \frac{\|u\|^2}{\dyle\left(\mu_3 \intn u^4dx\right)^{\frac12}}.
\ee
As a consequence, from \eqref{eq:inf1}, resp. \eqref{eq:inf2}, it follows that
\be\label{eq:pm1}
\frac{ \|(u_{1\e}^\pm,u_{2\e}^\pm)\|^2}{ \dyle\left(\mu_1 \intn (u_{1\e}^\pm)^4\,dx+\mu_2 \intn (u_{2\e}^\pm)^4\,dx+2\b\intn (u_{1\e}^\pm)^2(u_{2\e}^\pm)^2\,dx\right)^{\frac12}} \ge  \|(u^0,v^0)\|,
\ee
resp.
\be\label{eq:pm2}
\frac{\|u_{3\e}^\pm\|^2}{\dyle\left(\mu_3 \intn (u_{3\e}^\pm)^4dx\right)^{\frac 12}} \geq  \|U_3\|.
\ee
Multiplying  the third equation of $(S_3)$ by $u_{3\e}^\pm$ and integrating on $\Rn$ one infers
\begin{eqnarray*}
\|u_{3\e}^\pm\|^2 &=& \mu_3 \intn (u_{3\e}^\pm)^4dx + \e \intn\left[ (u_{3\e}^\pm)^2(\wt{\b}_{13} u_{1\e}^2 +\wt{\b}_{23} u_{2\e}^2)   \right]dx\\
&\leq&\mu_3 \intn (u_{3\e}^\pm)^4dx \\ & & +\e \left( \intn (u_{3\e}^\pm)^4dx\right)^{1/2}
\left[ \wt{\b}_{13}\left( \intn u_{1\e}^4\right)^{1/2} +\wt{\b}_{23}\left(\intn u_{2\e}^4\right)^{1/2}   \right].
\end{eqnarray*}
This, jointly with \eqref{eq:pm2}, yields
$$
\|u_{3\e}^\pm\|^2 \le \frac{\|u_{3\e}^\pm\|^4}{ \|U_3\|^2} + \e\, \vartheta_\e\;\frac{\|u_{3\e}^\pm\|^2}{ \|U_3\|},
$$
where
$$
\vartheta_\e = \mu_3^{-1/2}\left[\wt{\b}_{13}\left( \intn u_{1\e}^4\right)^{1/2} +\wt{\b}_{23}\left(\intn u_{2\e}^4\right)^{1/2}   \right].
$$
Since $\bu_\e\to \bu_0$, clearly  $(u_{1\e},u_{2\e})\to (u^0,v^0)$, then one has $\e\vartheta_\e\to 0$ as $\e\to 0$.
Hence, if  $\|u_{3\e}^\pm\|>0$, one obtains
\begin{equation}\label{eq:pmm1}
\|u_{3\e}^\pm\|^2 \ge  \|U_3\|^2 +o(1),
\end{equation}
where $o(1)=o_\e (1)\to 0$ as $\e\to 0$.
Using again  $\bu_\e \to \bu_0$, then $u_{3\e}\to U_3>0 $ and as a consequence, for $\e$ small enough, $\|u_{3\e}^+\|>0$. Thus
\eqref{eq:pmm1} gives
\begin{equation}\label{eq:pm11}
\|u_{3\e}^+\|^2 \ge  \|U_3\|^2 +o(1).
\end{equation}
Multiplying now the first, resp. the second equation of
$(S_3)$ by $u_{1\e}^\pm$, resp.  $u_{2\e}^\pm$ and integrating on $\Rn$ one infers
\begin{eqnarray*}
\|(u_{1\e}^\pm,u_{2\e}^\pm)\|^2 &=& \mu_1 \intn (u_{1\e}^\pm)^4dx + \mu_2 \intn (u_{2\e}^\pm)^4dx\\
& & +\b  \intn \left[(u_{1\e}^\pm)^2u_{2\e}^2+ (u_{2\e}^\pm)^2u_{1\e}^2\right]dx \\
 &  &+ \e \intn u_{3\e}^2\left[ (\wt{\b}_{13} (u_{1\e}^\pm)^2 +\wt{\b}_{23} (u_{2\e}^\pm)^2)   \right]dx\\
&\le &\mu_1 \intn (u_{1\e}^\pm)^4dx + \mu_2 \intn (u_{2\e}^\pm)^4dx+ 2\b  \intn (u_{1\e}^\pm)^2(u_{2\e}^\pm)^2 dx \\
& & +\b  \intn \left[(u_{1\e}^\mp)^2(u_{2\e}^\pm)^2+ (u_{1\e}^\pm)^2(u_{2\e}^\mp)^2\right]dx\\
& & +\e \left( \intn u_{3\e}^4\right)^{\frac 12}\left[ \wt{\b}_{13} \left(\intn(u_{1\e}^\pm)^4\right)^{\frac 12} +\wt{\b}_{23}\left(\intn (u_{2\e}^\pm)^4\right)^{\frac 12}   \right]\! dx.
\end{eqnarray*}
This, jointly with \eqref{eq:pm1}, yields
\be
\begin{array}{rcl}\label{eq:masmenos1}
\|(u_{1\e}^\pm,u_{2\e}^\pm)\|^2 &  \leq  & \dyle\frac{\|(u_{1\e}^\pm,u_{2\e}^\pm)\|^4}{ \|(u^0,v^0)\|^2} + \b  \intn \left[(u_{1\e}^\mp)^2(u_{2\e}^\pm)^2+ (u_{1\e}^\pm)^2(u_{2\e}^\mp)^2\right]dx\\ & & \\
& &\dyle +\e\, \psi_\e\;\frac{\|(u_{1\e}^\pm,u_{2\e}^\pm)\|^2}{ \|(u^0,v^0)\|},
\end{array}
\ee
where
$$
\psi_\e =C\left( \intn u_{3\e}^4\right)^{1/2} \quad\mbox{for some constant}\quad C=C(\m_1,\m_2,\l,\b)>0.
$$
Using  that $u_{3\e}\to U_3$ we obtain
\be\label{eq:masmenos2}
\e\psi_\e\to 0\quad\mbox{as }\e\to 0,
\ee
and, from $(u_{1\e},u_{2\e})\to (u^0,v^0)>\bo$, we get $(u_{1\e}^-,u_{2\e}^-)\to \bo$ which implies
\be\label{eq:masmenos3}
\b  \intn \left[(u_{1\e}^\mp)^2(u_{2\e}^\pm)^2+ (u_{1\e}^\pm)^2(u_{2\e}^\mp)^2\right]dx\to 0\quad \mbox{as }\e\to 0.
\ee
Since $(u_{1\e}^+, u_{2\e}^+)\to (u^0,v^0)>\bo $, we have $\|(u_{1\e}^+, u_{2\e}^+)\|>0$ for $\e$ sufficiently small, thus
\eqref{eq:masmenos1}-\eqref{eq:masmenos3} gives
\be\label{eq:pm-fin11}
\|(u_{1\e}^+ ,u_{2\e}^+)\|^2 \geq  \|(u^0,v^0)\|^2 +o(1).
\ee
If  $\|(u_{1\e}^-, u_{2\e}^-)\|>0$, it is not so easy to obtain a similar estimate like \eqref{eq:masmenos3} with $\|(u_{1\e}^+ ,u_{2\e}^+)\|^2$ replaced by $\|(u_{1\e}^- ,u_{2\e}^-)\|^2$. To do so, we need to estimate more carefully the mixed term by \eqref{eq:masmenos3} since it has the same order of decay as  $\|(u_{1\e}^-, u_{2\e}^-)\|^2$. 

By the Cauchy-Schwarz inequality, the convergence $\bu_\e\to \bu_0$ and Remark \ref{rem:sobolev-cte} we have
\be\begin{array}{rcl}\label{eq:new-estimate}
& &\dyle  \b \intn \left[(u_{1\e}^-)^2(u_{2\e}^+)^2+ (u_{1\e}^+)^2(u_{2\e}^-)^2\right]dx\\ & & \\
 & & \le \dyle
\b \left[\left( \intn (u_{1\e}^-)^4\right)^{1/2}\left( \intn (u_{2\e}^+)^4\right)^{1/2} + \left( \intn (u_{1\e}^+)^4\right)^{1/2}\left( \intn (u_{2\e}^-)^4\right)^{1/2}\right]\\ & & \\
& & \le  \dyle \b\frac{1}{\s_\l^2}\| u_{1\e}^{-}\|^2
\left(\l^{1-\frac{n}{4}}\left(\intn U^4\, dx\right)^\frac 12+o(1)\right)\frac{\b-\m_1}{\b^2-\m_1\m_2}
\\ & & \\
& &\quad \dyle + \b\frac{1}{\s_\l^2}\| u_{2\e}^{-}\|^2
\left(\l^{1-\frac{n}{4}}\left(\intn U^4\, dx\right)^\frac 12+o(1)\right)\frac{\b-\m_2}{\b^2-\m_1\m_2}\\
& & \\ & & \le  (\rho+o(1))\|(u_{1\e}^- ,u_{2\e}^-)\|^2,
\end{array}
\ee
where
\be\label{eq:rho}
\rho=\max\left\{ \frac{\b-\m_1}{\b^2-\m_1\m_2},\,\frac{\b-\m_2}{\b^2-\m_1\m_2}\right\}.
\ee
Even more, taking into account that $\b>\max\{\m_1,\m_2\}$ then $0<\rho<1$. As a consequence, if $\|(u_{1\e}^-, u_{2\e}^-)\|>0$, by  \eqref{eq:masmenos1}, \eqref{eq:masmenos2} and \eqref{eq:new-estimate}  we have the following estimate,
\be\label{eq:correjida}
\|(u_{1\e}^-, u_{2\e}^-)\|^2 \geq  (1-\rho)\|(u^0,v^0)\|^2 +o(1).
\ee
Now, suppose by contradiction that there exists $k\in \{1,2,3\}$ such that
$\|u_{k\e}^-\|>0$. Then we have two possibilities:
\begin{itemize}
\item
If $k=3$, as in \eqref{eq:pmm1} one obtains $\|u_{3\e}^-\|^2 \geq  \|U_3\|^2+o(1)$, hence
\begin{equation}\label{eq:pm12}
\|\bu_\e^-\|^2 =\|(u_{1\e}^- ,u_{2\e}^-)\|^2+\|u_{3\e}^-\|^2 \geq  \|U_3\|^2+o(1).
\end{equation}
Next, we evaluate the functional
\be\label{eq:eval-funct}
\Phi(\bu_\e)= \tfrac 14 \|\bu_\e\|^2 = \tfrac 14 \left[  \|\bu_\e^+\|^2 + \|\bu_\e^-\|^2\right].
\ee
On one hand, using \eqref{eq:pm11},  \eqref{eq:pm-fin11}  and \eqref{eq:pm12}, we infer
\begin{equation}\label{eq:pmm2}
\Phi(\bu_\e)\geq \tfrac 14 \|\bu_0\|^2 +\tfrac 14\|U_3\|^2+o(1).
\end{equation}
On the other hand, since $\bu_\e \to \bu_0$ we also find
\be\label{eq:final}
\Phi(\bu_\e)= \tfrac 14 \|\bu_\e\|^2\to \tfrac 14\|\bu_0\|^2,
\ee
which is in contradiction with \eqref{eq:pmm2}, proving that $u_{3\e}\geq 0$.
\item If $k\in \{1,2\}$, by \eqref{eq:correjida} we have
\begin{equation}\label{eq:pm122}
\|\bu_\e^-\|^2 =\|(u_{1\e}^- ,u_{2\e}^-)\|^2+\|u_{3\e}^-\|^2 \ge (1-\rho)\|(u^0,v^0)\|^2+o(1).
\end{equation}
Using \eqref{eq:pm-fin11}, \eqref{eq:correjida}, \eqref{eq:eval-funct} and \eqref{eq:pm122}, we get
$$
\Phi(\bu_\e)\geq \tfrac 14 \|\bu_0\|^2 +\tfrac 14(1-\rho)\|(u^0,v^0)\|^2+o(1).
$$
This is a contradiction with \eqref{eq:final}, proving that $(u_{1\e},u_{2\e})\ge \bo$.
\end{itemize}
In conclusion, we have proved that $\bu_\e\ge \bo$. Finally, if $\b_{12},\, \b_{23}\ge 0$, using once more that $\bu_\e\to \bu_0$ and applying the maximum principle it follows that  $\bu_\e>\bo$.
\rule{2mm}{2mm}

\

Let us set
\be\label{eq:limite}
\bz=((u_1^0,v_1^0),\ldots,(u_m^0,v_m^0),U_{2m+1},\ldots,U_{N}),
\ee
where $(u_k^0,v_k^0)$ is given by \eqref{eq:2-sistema} with the parameters $\l=\l_{2k}$, $\m_1=\m_{2k-1}$, $\m_2=\m_{2k}$, $\b=\b_{k}$, and $1\le k\le m$; see Theorem A.
And
$$
U_{2m+\ell}(x)=\sqrt{\dfrac{\l_{2m+\ell}}{\m_{2m+\ell}}}U(\sqrt{\l_{2m+\ell}}x)
$$
where $U$ is given by \eqref{eq:single}; see Remark \ref{rem:1}.

The main result of the paper is the following.
\begin{Theorem}\label{th:ec2}
If $\b_k>\max \{ \m_{2k-1},\m_{2k}\}$ for  $1\le k\le m$, there exists $\e_0>0$ such that for $0<\e<\e_0$, \eqref{eq:main-sys} has a radial bound state $\bu_{\e}\ge \bo$  with $\bu_{\e}\to \bz$ as $\e\to 0$. Furthermore, if  $\b_{m+l}\ge 0$ for all $1\le\ell\le N-2m$ then $\bu_{\e}>\bo$.
\end{Theorem}
\begin{Remark}
Notice that a condition on the coupling factors like in Remark \ref{rem:2} holds provided $\e_0$  is small enough.
\end{Remark}
\noindent {\it Proof of Theorem \ref{th:ec2}}.
We follow, with suitable modifications, the arguments in the proof of Theorem \ref{th:ec1}.

Let us consider the critical point $\bz$  (defined in \eqref{eq:limite}) of the unperturbed functional $\Phi_0$. We note that $U_{2m+\ell}$ is a non-degenerate critical point of $I_{2m+\ell}$ on $H$ for each $\ell$; see \cite{kw} and Remark \ref{rem:1}. Also, $(u_k^0,v_k^0)$ is a non-degenerate critical point of the corresponding functional $\Phi_{k}$ described in \eqref{eq:Phi-k}, for each $1\le k\le m$; see Remark \ref{rem:non-degenerate}.
As an immediate consequence, $\bu_0$ is a non-degenerate critical point of $\Phi_0$ on $\mathbb{H}$, then  the Local Inversion Theorem provides us with the existence of a critical point $\bu_\e$ of $\Phi_\e$ for any $\e<\e_0$ with $\e_0$ sufficiently small. Furthermore, $\bu_\e\to \bu_0$ on $\mathbb{H}$ as $\e\to 0$. To complete the proof it remains to show that $\bz\ge \bo$ and if  $\b_{m+l}\ge 0$ for all $1\le\ell\le N-2m$, then $\bu_{\e}>\bo$.

Following the notation of $\bz$ by \eqref{eq:limite}, let us denote
$$
\bu_\e=((u_1^\e,v_1^\e),\ldots,(u_m^\e,v_m^\e), u^\e_{2m+1},\ldots,u_N^\e).
$$
Using that $\bu_\e\to \bz>\bo$, as in \eqref{eq:pm11} and \eqref{eq:pm-fin11} it follows
\be\label{eq:pm-final1}
\|((u_{k}^\e)^+,(v_{k}^\e)^+)\|^2\ge \|(u_{k}^0,v_{k}^0)\|^2+o(1)\quad \mbox{for any }1\le k\le m,
\ee
 and
\be\label{eq:pm-final2}
 \|(u_{2m+\ell}^\e)^+\|^2\ge \|U_{2m+\ell}\|^2+o(1) \quad \mbox{for any }1\le\ell\le N-2m.
\ee
Suppose  by contradiction that there exists either $1\le k_0\le m$, or $1\le \ell_0\le N-2m$ such that either
\be\label{eq:either}
\|((u_{k_0}^\e)^-,(v_{k_0}^\e)^-)\|>0,
\ee
or
\be\label{eq:or}
\|(u_{2m+\ell_0}^\e)^-\|>0.
\ee
\begin{itemize}
\item If \eqref{eq:either} holds, by \eqref{eq:correjida}
\be\label{eq:pm-final1menos}
\|((u_{k_0}^\e)^- ,(v_{k_0}^\e)^-)\|^2 \ge  (1-\rho_{k_0})\|(u_{k_0}^0,v_{k_0}^0)\|^2 +o(1),
\ee
for some $0<\rho_{k_0}<1$ defined by \eqref{eq:rho} with with the parameters $\l=\l_{2{k_0}}$, $\m_1=\m_{2{k_0}-1}$, $\m_2=\m_{2k_0}$, $\b=\b_{k_0}$. Then we obtain
\be\label{eq:either2}
\|\bu_\e^-\|^2  \ge (1-\rho_{k_0})\|(u_{k_0}^0,v_{k_0}^0)\|^2 +o(1),
\ee
and by  \eqref{eq:pm-final1}, \eqref{eq:pm-final1menos} we have the following inequality
\be\label{eq:either22}
\Phi(\bu_\e)\ge \tfrac 14 \|\bz\|^2 +\tfrac 14 (1-\rho_{k_0}) \|(u_{k_0}^0,v_{k_0}^0)\|^2 +o(1).
\ee
\item If \eqref{eq:or} holds, by \eqref{eq:pm11}
\be\label{eq:pm-final2menos}
\|(u_{2m+\ell_0}^\e)^-\|^2 \geq  \|U_{2m+\ell_0}\|^2 +o(1),
\ee
which implies
\be\label{eq:or2}
\|\bu_\e^-\|^2  \ge \|U_{2m+\ell_0}\|^2 +o(1),
\ee
and by \eqref{eq:pm-final2} and \eqref{eq:pm-final2menos} the following inequality holds
\be\label{eq:or22}
\Phi(\bu_\e)\ge \tfrac 14 \|\bz\|^2 +\tfrac 14 \|U_{2m+\ell_0}\|^2 +o(1).
\ee
\end{itemize}
Since $\bu_\e\to \bz$ we have
$$
\Phi(\bu_\e)=\tfrac 14 \|\bu_\e\|^2\to \tfrac 14 \|\bz\|^2.
$$
This is a contradiction with \eqref{eq:either22} and also with \eqref{eq:or22}, proving that $\bu_\e\ge \bo$. To finish,  if  $\b_{m+l}\ge 0$ for all $1\le\ell\le N-2m$,  using once more that $\bu_\e\to \bz$ and applying the maximum principle it follows that  $\bu_\e>\bo$.
\rule{2mm}{2mm}

\

If $2m=N$, or equivalently $\ell=0$, we set $\bz_1=((u_1^0,v_1^0),\ldots,(u_m^0,v_m^0))$, and if $m=0$, or equivalently $\ell=N$, we set $\bz_2=(U_1,\ldots,U_N)$. Then as a consequence of Theorem \ref{th:ec2} we have the following.
\begin{corollary}\label{cor:ec1}
For $\e$  small enough,
\begin{itemize}
\item[(i)] if $2m=N\Leftrightarrow\ell=0$, then  \eqref{eq:main-sys} has a radial bound state $\bu_{\e}^{(1)}>\bo$  and  $\bu_{\e}^{(1)}\to \bz_1$ as $\e\to 0$.
\item[(ii)] if $m=0\Leftrightarrow \ell=N$,  then \eqref{eq:main-sys} has a radial bound state $\bu_{\e}^{(2)}\ge \bo$  and $\bu_{\e}^{(2)}\to \bz_2$ as $\e\to 0$. Furthermore, if  $\b_{\ell}\ge 0$ for all $1\le\ell\le N$ then $\bu_{\e}^{(2)}>\bo$.
\end{itemize}
\end{corollary}
We observe that Corollary \ref{cor:ec1}-(ii)  is nothing but \cite[Theorem 6.4]{ac2}.
\begin{Remark}
Note that in Theorem \ref{th:ec2}  some of the coupling factors $\b_k$, $1\le k\le m$ are big while the others $\b_{k+\ell}$, $1\le \ell \le N-2m$ are small, and  some of these can be negative giving rise to a non-negative new non-trivial bound states. Similar remark can be done for Corollary \ref{cor:ec1}-(ii).
\end{Remark}

\

\noindent {\bf Acknowledgements.}
The author thought about this work during a visit  at the Institut f\"ur Mathematik of the Goethe-Universit\"at Frankfurt in 2013, and he wishes to thank its hospitality. That visit was partially supported by Beca de Movilidad de Profesores de las Universidades P\'ublicas de Madrid 2012-2013 of Fundaci\'on Caja Madrid. The author also wants to thank Prof. Tobias Weth for useful discussions during that visit.


{\small
\begin{thebibliography}{777}
\bibitem{akanbook} N. Akhmediev, A. Ankiewicz,  ``Solitons,
Nonlinear pulses and beams''. Champman \& Hall, London, 1997.

\bibitem{acv} P. \'Alvarez-Caudevilla, E. Colorado, V. A. Galaktionov, {\it Existence of solutions for a system of coupled nonlinear stationary bi-harmonic Schr\"odinger equations}, Preprint arXiv:1402.4165v2, 2014.

\bibitem{ac1} A. Ambrosetti, E. Colorado, {\it Bound and ground states of coupled nonlinear Schr\"odinger equations.} C. R. Math.
Acad. Sci. Paris {\bf 342} (2006), no. 7, 453-458.

\bibitem{ac2}  A. Ambrosetti, E. Colorado, {\it Standing waves of some coupled nonlinear Schr\"odin\-ger equations.} J. Lond. Math.
Soc. ({\bf 2}) 75 (2007), no. 1, 67-82.

\bibitem{acr} A. Ambrosetti, E. Colorado, D. Ruiz, {\it Multi-bump solitons to linearly coupled systems of nonlinear Schr\"odinger
equations.} Calc. Var. Partial Differential Equations {\bf 30} (2007), no. 1, 85-112.

\bibitem{a-m} A. Ambrosetti and A. Malchiodi,  ``Perturbation methods and semilinear elliptic problems on $\Rn$", Progress in Math. Vol. {\bf 240}, Birkh\"auser, 2005.

\bibitem{bt} T. Bartsch, Z.-Q. Wang, {\it Note on ground states of nonlinear Schr\"odinger systems.} J. Partial Differential
Equations {\bf 19} (2006), no. 3, 200-207.

\bibitem{chen-zou} Z. Chen; W. Zou, {\it An optimal constant for the existence of least energy solutions of a coupled Schrödinger system.} Calc. Var. Partial Differential Equations {\bf 48} (2013), no. 3-4, 695–711.

\bibitem{cing} S. Cingolani, {\it Positive solutions to perturbed elliptic problems in $\Rn$ involving critical Sobolev exponent}, Nonlin. Anal. T.M.A. {\bf 48}  (2002), 1165-1178.

\bibitem{col} E. Colorado, {\it Existence results for some systems of coupled Fractional Nonlinear Schr\"odinger Equations.} Recent trends in nonlinear partial differential equations. II. Stationary problems, 135-150, Contemp. Math., {\bf 595}, Amer. Math. Soc., Providence, RI, 2013.

\bibitem{dancer-wei} E. N.  Dancer, J. Wei, {\it Spike solutions in coupled nonlinear Schr\"odinger equations with attractive interaction}. Trans. Amer. Math. Soc. {\bf 361} (2009), no. 3, 1189-1208.

\bibitem{hioe} F.T. Hioe, T.S. Salter, {\it Special set and solutions of coupled nonlinear
Schr\"odinger equations}. J. Phys. A: Math. Gen. {\bf 35} (2002), 8913-8928.

\bibitem{itanaka} N. Ikoma, K.  Tanaka, {\it A local mountain pass type result for a system of nonlinear Schr\"odinger equations.}
Calc. Var. Partial Differential Equations {\bf 40} (2011), no. 3-4, 449-480.

\bibitem{kw} M.K. Kwong, {\it Uniqueness of positive solutions of $\Delta u -u + u^p = 0$ in $\R^N$}. Arch. Rat. Mech. Anal.
\textbf{105} (1989), 243-266.

\bibitem{linwei} T-C. Lin, J. Wei,  {\it Ground state of $N$ coupled nonlinear Schr\"odinger  equations in $\R^n$, $n\leq 3$}.
Comm. Math. Phys. {\bf 255} (2005), 629-653.

\bibitem{lions} P.L. Lions, {\it  Sym\'etrie et compacit\'e dans les espaces de Sobolev.} J. Funct. Anal. {\bf 49} (1982), no. 3, 315-334.

\bibitem{liu-wang} Z. Liu, Z.-Q. Wang, {\it Ground states and bound states of a nonlinear Schrödinger system.} Adv. Nonlinear Stud. {\bf 10} (2010), no. 1, 175-193.

\bibitem{mmp} L. Maia, E. Montefusco, B. Pellacci, {\it  Positive solutions for a weakly coupled nonlinear Schr\"odinger system.}
J. Differential Equations {\bf 229} (2006), no. 2, 743-767.

\bibitem{me} C. R. Menyuk, {\it Nonlinear pulse propagation in birifrangent optical fibers}. IEEE Jour. Quantum Electr.
{\bf 23}-2 (1987), 174-176.

\bibitem{pomp} A. Pomponio, {\it Coupled nonlinear Schr\"odinger systems with potentials}. J. Differential Equations
{\bf 227} (2006), no. 1, 258-281.

\bibitem{sirakov} B. Sirakov, {\it  Least energy solitary waves for a system of nonlinear Schr\"odinger equations in
$\Bbb R^n$.} Comm. Math. Phys. {\bf 271} (2007), no. 1, 199-221.

\bibitem{soave} N. Soave, {\it On existence and phase separation of solitary waves for nonlinear Schrödinger systems modelling simultaneous cooperation and competition}.
Preprint,  arXiv:1310.8492v2, 2013.

\bibitem{wei-yao} J. Wei, W. Yao {\it Uniqueness of positive solutions to some coupled nonlinear Schr\"odinger
equations.} Commun. Pure Appl. Anal. {\bf 11} (2012), no. 3, 1003-1011.

\bibitem{ww} J. Wei, T. Weth {\it Nonradial symmetric bound states for a system of coupled Schrödinger equations.}
Atti Accad. Naz. Lincei Cl. Sci. Fis. Mat. Natur. Rend. Lincei (9) Mat. Appl. {\bf 18} (2007), no. 3, 279-293.
\end{thebibliography}}
\enddocument